\documentclass{article}
\usepackage{amsmath,amsthm,amssymb,amsfonts}
\usepackage{verbatim}
\usepackage{graphicx}
\usepackage{stmaryrd} 
\usepackage{hyperref}

\usepackage[colorinlistoftodos]{todonotes}

\theoremstyle{plain}
\newtheorem{thm}{Theorem}

\newtheorem{prop}[thm]{Proposition}
\newtheorem{cor}[thm]{Corollary}

\newtheorem{remark}[thm]{Remark}
\newtheorem{question}[thm]{Question}

\theoremstyle{definition}
\newtheorem{definition}[thm]{Definition}
\newtheorem{exl}[thm]{Example}

\numberwithin{thm}{section}

\newcommand{\adj}{\leftrightarrow}
\newcommand{\adjeq}{\leftrightarroweq}

\DeclareMathOperator{\id}{id}
\def\Z{{\mathbb Z}}

\def\R{{\mathbb R}}

\begin{document}
\title{Convexity and AFPP in the Digital Plane}
\author{Laurence Boxer
\thanks{
    Department of Computer and Information Sciences,
    Niagara University,
    Niagara University, NY 14109, USA;
    and Department of Computer Science and Engineering,
    State University of New York at Buffalo.
    email: boxer@niagara.edu
}
}

\date{ }
\maketitle{}

\begin{abstract}
We examine the relationship between convexity and the approximate fixed point property (AFPP)
for digital images in $\Z^2$.

\end{abstract}


\section{Introduction}
The study of fixed points is prominent in many branches of mathematics. In digital topology,
it has become worthwhile to broaden the study to ``approximate fixed points." The Approximate
Fixed Point Property (AFPP), a generalization of the classical fixed point property (FPP),
was introduced in~\cite{BEKLL}. In this paper, we show that
for digital images $X \subset \Z^2$, convexity can help us show whether $(X,c_2)$ has
the AFPP.

\section{Preliminaries}
Much of this section is quoted or paraphrased from papers that are listed in the
references, especially~\cite{BxAFPP,BxAFPPtreesProducts,BxConvexity,BEKLL}.

We use $\Z$ to indicate the set of integers; $\R$ for the set of real numbers.

For $(x,y) \in \Z^2$, the projection functions $pr_1,pr_2: \Z^2 \to \Z$ are
\[ pr_1(x,y) = x, ~~~ pr_2(x,y) = y.
\]

\subsection{Adjacencies}
A digital image is a graph $(X,\kappa)$, where $X$ is a subset of $\Z^n$ for
some positive integer~$n$, and $\kappa$ is an adjacency relation for the points
of~$X$. The $c_u$-adjacencies are commonly used.
Let $x,y \in \Z^n$, $x \neq y$, where we consider these points as $n$-tuples of integers:
\[ x=(x_1,\ldots, x_n),~~~y=(y_1,\ldots,y_n).
\]
Let $u \in \Z$,
$1 \leq u \leq n$. We say $x$ and $y$ are 
{\em $c_u$-adjacent} if
\begin{itemize}
\item There are at most $u$ indices $i$ for which 
      $|x_i - y_i| = 1$.
\item For all indices $j$ such that $|x_j - y_j| \neq 1$ we
      have $x_j=y_j$.
\end{itemize}
Often, a $c_u$-adjacency is denoted by the number of points
adjacent to a given point in $\Z^n$ using this adjacency.
E.g.,
\begin{itemize}
\item In $\Z^1$, $c_1$-adjacency is 2-adjacency.
\item In $\Z^2$, $c_1$-adjacency is 4-adjacency and
      $c_2$-adjacency is 8-adjacency.
\item In $\Z^3$, $c_1$-adjacency is 6-adjacency,
      $c_2$-adjacency is 18-adjacency, and $c_3$-adjacency
      is 26-adjacency.
\end{itemize}

We write $x \adj_{\kappa} x'$, or $x \adj x'$ when $\kappa$ is understood, to indicate
that $x$ and $x'$ are $\kappa$-adjacent. Similarly, we
write $x \adjeq_{\kappa} x'$, or $x \adjeq x'$ when $\kappa$ is understood, to indicate
that $x$ and $x'$ are $\kappa$-adjacent or equal.

A subset $Y$ of a digital image $(X,\kappa)$ is
{\em $\kappa$-connected}~\cite{Rosenfeld},
or {\em connected} when $\kappa$
is understood, if for every pair of points $a,b \in Y$ there
exists a sequence $\{y_i\}_{i=0}^m \subset Y$ such that
$a=y_0$, $b=y_m$, and $y_i \adj_{\kappa} y_{i+1}$ for $0 \leq i < m$.

Given a digital image $(X,\kappa)$ and $x \in X$, we denote by $N^*(X,\kappa,x)$ the set
$\{y \in X \, | \, y \adjeq_{\kappa} x\}$.

\subsection{Digitally continuous functions}
The following generalizes a definition of~\cite{Rosenfeld}.

\begin{definition}
\label{continuous}
{\rm ~\cite{Boxer99}}
Let $(X,\kappa)$ and $(Y,\lambda)$ be digital images. A single-valued function
$f: X \rightarrow Y$ is $(\kappa,\lambda)$-continuous if for
every $\kappa$-connected $A \subset X$ we have that
$f(A)$ is a $\lambda$-connected subset of $Y$. $\Box$
\end{definition}

When the adjacency relations are understood, we will simply say that $f$ is \emph{continuous}. Continuity can be expressed in terms of adjacency of points:
\begin{thm}
{\rm ~\cite{Rosenfeld,Boxer99}}
A function $f:X\to Y$ is continuous if and only if $x \adj x'$ in $X$ implies 
$f(x) \adjeq f(x')$. \qed
\end{thm}

See also~\cite{Chen94,Chen04}, where similar notions are referred to as {\em immersions}, {\em gradually varied operators},
and {\em gradually varied mappings}.

Let $Y \subset X$ and let $f: X \to Y$ be $(\kappa,\kappa)$-continuous such that
$r(y)=y$ for all $y \in Y$. Then $r$ is a {\em $\kappa$-retraction}.

The notation $C(X,\kappa)$ denotes
$\{f: X \to X \, | \, f \mbox{ is }(\kappa,\kappa)-\mbox{continuous}\}$.

\subsection{Approximate fixed points and the AFPP}
\label{approxPrelim}
Let $f \in C(X,\kappa)$
and let $x \in X$. We say
\begin{itemize}
    \item $x$ is a {\em fixed point} of $f$ if $f(x)=x$; 
    \item If $f(x) \adjeq_{\kappa} x$, then
          $x$ is an {\em almost fixed point}~\cite{Rosenfeld,Tsaur} or
          {\em approximate fixed point}~\cite{BEKLL} of 
          $(f,\kappa)$.
    \item A digital image $(X,\kappa)$ has the
          {\em approximate fixed point property} (AFPP)~\cite{BEKLL} if for every $f \in C(X,\kappa)$
          there is an approximate fixed point of $f$. This generalizes the {\em fixed point property}
          (FPP): a digital image $(X,\kappa)$ has the FPP if every $f \in C(X,\kappa)$ has a
          fixed point.
\end{itemize}

The AFPP gathered attention in part because only a digital image with a single point has the
FPP~\cite{BEKLL}.

A. Rosenfeld's paper~\cite{Rosenfeld} states the following as its Theorem~4.1 (quoted verbatim).
\begin{quote}
    Let $I$ be a digital picture, and let $f$ be a continuous function from $I$
    into $I$; then there exists a point $P \in I$ such that $f(P)=P$ or is a neighbor
    or diagonal neighbor of $P$.
\end{quote}
We quote from~\cite{BxAFPP}:
\begin{quote}
Several subsequent papers have incorrectly
concluded that this [Rosenfeld's] result implies that $I$ with
some $c_u$ adjacency has the $AFPP_S$. 
By {\em digital picture} Rosenfeld means a digital cube, $I= [0,n]_{\Z}^v$.
By a ``continuous function" he means a $(c_1,c_1)$-continuous function;
by ``a neighbor or diagonal neighbor of $P$" he means a $c_v$-adjacent point.
\end{quote}
Thus, Rosenfeld's result was important but weaker than that of Theorem~\ref{haveAFPP}(6), below.

\begin{thm}
\label{haveAFPP}
The following digital images have the AFPP.
\begin{enumerate}
    \item Any digital interval $([a,b]_{\Z}, c_1)$~\cite{Rosenfeld,BEKLL}.
    \item Any digital image $(Y,\lambda)$ that is isomorphic to
          $(X,\kappa)$ such that $(X,\kappa)$ has the AFPP~\cite{BEKLL}
    \item Any digital image $(Y,\kappa)$ that is a retract of
          $(X,\kappa)$ such that $(X,\kappa)$ has the AFPP~\cite{BEKLL}.
    \item Any digital image $(T,\kappa)$ that is a tree~\cite{BxAFPPtreesProducts}.
    \item Any digital image $(X,c_{m+n})$ such that $X = X' \times \Pi_{i=1}^n [a_i,b_i]_{\Z}$,
           $X' \subset \Z^m$, and $(X',c_m)$ has the AFPP~\cite{BxAFPPtreesProducts}.
    \item Any digital cube $(\Pi_{i=1}^n [a_i,b_i]_{\Z}, c_n)$~\cite{BxAFPPtreesProducts}.
\end{enumerate}
\end{thm}

The next result suggests that ``most" digital images $(X,c_u) \subset \Z^v$ that have the
AFPP have $u=v$.

\begin{thm}
{\rm \cite{BxAFPP}}
Let $X \subset \Z^v$ be such that $X$ has a subset $Y = \Pi_{i=1}^v [a_i,b_i]_{\Z}$, where
$v>1$; for all indices $i$, $b_i \in \{a_i,a_i + 1\}$; and, for at least 2 indices $i$,
$b_i=a_i+1$. Then $(X,c_u)$ fails to have the AFPP for $1 \le u < v$.
\end{thm}

\begin{exl}
{\rm \cite{BEKLL}}
\label{SCC4notAFPP}
A digital simple closed curve of at least 4 points does not have the AFPP.
\end{exl}

\subsection{Digital convexity, disks}
Material in this section is quoted or paraphrased from~\cite{BxConvexity}.

Let $\kappa \in \{c_1,c_2\}$, $n>1$. We say a $\kappa$-connected 
set $S=\{x_i\}_{i=1}^n \subset \Z^2$ is a
{\em (digital) line segment} if the members of $S$ are collinear.

\begin{remark}
\label{segSlope}
{\rm \cite{BxConvexity}}
A digital line segment must be vertical, horizontal, or have
slope of $\pm 1$. We say a segment with slope of $\pm 1$ is
{\em slanted}.
\end{remark}

A {\em (digital) $\kappa$-closed curve} is a
path $S=\{s_i\}_{i=0}^{m-1}$ such that $i \neq j$ implies $s_i \neq s_j$,
and $s_i \adj_{\kappa} s_{(i+1) \mod \, m}$ for $0 \leq i \leq m-1$,
If $s_i \adj_{\kappa} s_j$ implies 
$j = (i \pm 1) \mod m$, $S$ is a {\em (digital) 
$\kappa$-simple closed curve}.
For a simple closed curve $S \subset \Z^2$ we generally assume
\begin{itemize}
    \item $m \ge 8$ if $\kappa = c_1$, and
    \item $m \ge 4$ if $\kappa = c_2$.
\end{itemize}
These requirements are necessary for the Jordan Curve
Theorem of digital topology, below, as a
$c_1$-simple closed curve in $\Z^2$ needs at least 8 points to
have a nonempty finite complementary $c_2$-component,
and a $c_2$-simple closed curve in $\Z^2$ needs at least 4 points to
have a nonempty finite complementary $c_1$-component.
Examples in~\cite{RosenfeldMAA} show why it is
desirable to consider $S$ and $\Z^2 \setminus S$
with different adjacencies.

\begin{thm}
{\rm \cite{RosenfeldMAA}}
{\em (Jordan Curve Theorem for digital topology)}
Let $\{\kappa, \kappa'\} = \{c_1, c_2\}$.
Let $S \subset \Z^2$ be a simple closed 
$\kappa$-curve such that $S$ has at least 8 points if
$\kappa = c_1$ and such that $S$ has at least 
4 points if $\kappa = c_2$. Then
$\Z^2 \setminus S$ has exactly 2 $\kappa'$-connected
components.
\end{thm}

One of the $\kappa'$-components of 
$\Z^2 \setminus S$ is finite and the other is infinite. This 
suggests the following.
\begin{definition}
\label{diskDef}
{\rm \cite{BxConvexity}}
Let $S \subset \Z^2$ be a $c_2$-closed curve such that
$\Z^2 \setminus S$ has two $c_1$-components, one finite and the
other infinite. The union $D$ of $S$ and the finite $c_1$-component 
of $\Z^2 \setminus S$ is a {\em (digital) disk}. $S$ is
a {\em bounding curve} of $D$. The finite $c_1$-component 
of $\Z^2 \setminus S$ is the {\em interior of} $S$, denoted $Int(S)$,
and the infinite $c_1$-component of $\Z^2 \setminus S$ is the {\em exterior of} 
$S$, denoted $Ext(S)$.
\end{definition}

Note a disk may have multiple distinct bounding curves~\cite{BxConvexity}.

More generally, we have the following.
\begin{definition}
\label{boundingCurvesDef}
{\rm \cite{BxConvexity}}
let $X \subset \Z^2$ be a finite, $c_i$-connected set,
$i \in \{1,2\}$. Suppose there are 
pairwise disjoint $c_2$-closed curves
$S_j \subset X$, $1 \le j \le n$, such that
\begin{itemize}
    \item $X \subset S_1 \cup Int(S_1)$;
    \item for $j>1$, $D_j = S_j \cup Int(S_j)$ is a digital disk;
    \item no two of 
    \[ S_1 \cup Ext(S_1), D_2, \ldots, D_n
    \]
    are $c_1$-adjacent or $c_2$-adjacent; and
    \item we have 
    \[ \Z^2 \setminus X = Ext(S_1) \cup \bigcup_{j=2}^n Int(S_j).
    \]
\end{itemize}
Then $\{S_j\}_{j=1}^n$ is a {\em set of bounding curves of} $X$.
\end{definition}

As above, $X$ may have multiple distinct sets of bounding curves.

A set $X$ in a Euclidean space $\R^n$ is
{\em convex} if for every pair of distinct
points $x,y \in X$, the line segment
$\overline{xy}$ from $x$ to $y$ is contained in $X$.
The {\em convex hull of} $Y \subset \R^n$,
denoted $hull(Y)$, is the
smallest convex subset of $\R^n$ that contains~$Y$.
If $Y \subset \R^2$ is a finite set, then
$hull(Y)$ is a single point if $Y$ is a singleton;
a line segment if $Y$ has at least 2 members and all are
collinear; otherwise, $hull(Y)$ is a polygonal disk,
and the endpoints of the edges of $hull(Y)$ are its {\em vertices}.

A digital version of convexity can be stated
for subsets of the digital plane~$\Z^2$ as follows.
A finite set $Y \subset \Z^2$ is 
{\em (digitally) convex} \cite{BxConvexity} if either
\begin{itemize}
    \item $Y$ is a single point, or
    \item $Y$ is a digital line segment, or
    \item $Y$ is a digital disk with a bounding curve $S$
          such that the endpoints of the maximal digital line segments
          of~$S$ are the vertices of $hull(Y) \subset \R^2$.
\end{itemize}

\section{Retractions, convexity, and the AFPP}
Due to assertions~(3) and~(6) of Theorem~\ref{haveAFPP}, the following theorem can
be useful in determining whether $(X,c_2)$ has the AFPP, for $X \subset \Z^2$.

\begin{thm}
\label{cvxRetractOfRect}
Let $X \subset Y=[a,b]_{\Z} \times [c,d]_{\Z} \subset \Z^2$, such that $X$ is a
digitally convex disk. Let $S$ be a bounding curve for $X$.
Then there is a $c_2$-retraction $r: Y \to X$ such that $r(Y \setminus Int(S)) = S$.
\end{thm}

\begin{figure}
    \centering
    \includegraphics[height=3.5in]{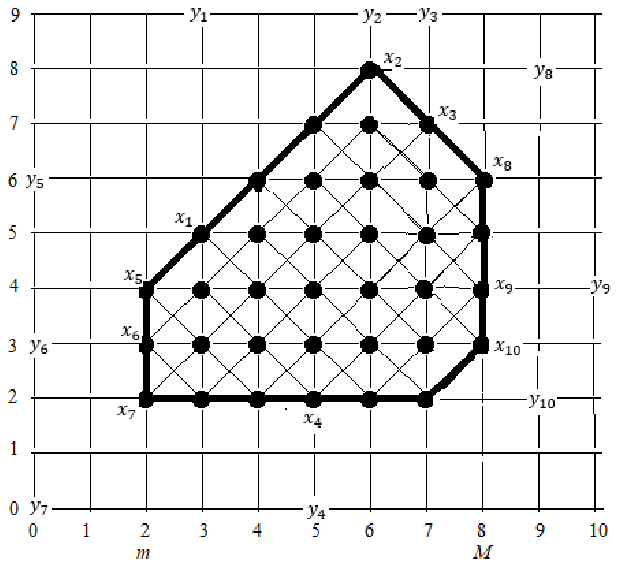}
    \caption{Retraction $r$ of a digital image $Y$ to a subset $X$ that is a convex disk as in
    Theorem~\ref{cvxRetractOfRect}. Here, $s_0=2$, $s_1=4$, $s_2=3$, $s_3=6$. \newline
    a) Each point vertically above or below the disk is mapped to its nearest 
       vertical neighbor in $X$, e.g., $r(y_i)=x_i$, $i \in \{1,2,3,4\}$. \newline
    b) Each point to the left (not necessarily horizontally) of $X$ is mapped to 
       the nearest member of $X$ with minimal first coordinate,
       e.g., $r(y_i)=x_i$, $i \in \{5,6,7\}$. \newline
    c) Each point to the right (not necessarily horizontally) of $X$ is mapped to 
       the nearest member of $X$ with maximal first coordinate,
       e.g., $r(y_i)=x_i$, $i \in \{8,9,10\}$.
       }
    \label{fig:convexRetract}
\end{figure}

\begin{proof}
We define a function $r: Y \to X$ as follows. For $y \in X$, $r(y)=y$.

For $y \not \in X$ we proceed as follows. Let
\[ m = \min\{p_1(x) \, | \, x \in X\}, ~~~ L = \{ x \in S \, | \, p_1(x) = m\},
\]
\[ M = \max\{p_1(x) \, | \, x \in X\}, ~~~ R = \{ x \in S \, | \, p_1(x) = R\}.
\]
\begin{itemize}
    \item Suppose $y=(a,b)$ for $m \le a \le M$. 
          \begin{itemize} 
             \item If $x = (a,c) \in S$ such that 
                   $b > c = \max\{n \, | \, (a,n) \in X\}$ then $r(y)=x$. 
                   (See $y_1,y_2,y_3$ in Figure~\ref{fig:convexRetract}.)
             \item If $x = (a,d) \in S$ such that 
                   $b < d = \min\{n \, | \, (a,n) \in X\}$ then $r(y)=x$.
                    (See $y_4$ in Figure~\ref{fig:convexRetract}.)
          \end{itemize}
          \item Suppose $y=(a,b)$ for $a < m$; then there is a unique
                nearest (in the Euclidean metric) $y' \in L$ to $y$, determined as follows. Let 
                \[ s_0 = \min\{pr_2(x) \, | \, x \in L\}, ~~~ s_1 = \max\{pr_2(x) \, | \, x \in L\}.
                \]
                \begin{itemize}
                    \item If $b > s_1$ then $y' = (m,s_1)$. (See $y_5$ in Figure~\ref{fig:convexRetract}.)
                    \item If $s_0 \le b \le s_1$ then $y'=(m,b)$. (See $y_6$ in Figure~\ref{fig:convexRetract}.)
                    \item If $b < s_0$ then $y'=(m,s_0)$. (See $y_7$ in Figure~\ref{fig:convexRetract}.)
                \end{itemize}
                Let $r(y)=y'$.
          \item Suppose $y=(a,b)$ for $a > M$; then there is a unique
                nearest (in the Euclidean metric) $y' \in L$ to $y$, determined as follows. Let 
                \[ s_2 = \min\{pr_2(x) \, | \, x \in R\}, ~~~ s_3 = \max\{pr_2(x) \, | \, x \in R\}.
                \]
                \begin{itemize}
                 \item If $b > s_3$ then $y' = (M,s_3)$. (See $y_8$ in Figure~\ref{fig:convexRetract}.)
                 \item If $s_2 \le b \le s_3$ then $y'=(M,b)$.  (See $y_9$ in Figure~\ref{fig:convexRetract}.)
                 \item If $b < s_2$ then $y'=(M,s_2)$.  (See $y_{10}$ in Figure~\ref{fig:convexRetract}.)
                \end{itemize}
                Let $r(y)=y'$.
\end{itemize}

In order to show $r$ is a $c_2$-retraction, we must show $r \in C(X,c_2)$. Let
$y \adj_{c_2} y'$ in $X$.
\begin{itemize}
    \item Suppose $y \in X$. 
          \begin{itemize}
              \item If $y' \in X$, then we have $r(y') = y' \adj_{c_2} y = r(y)$.
              \item If $y' \not \in X$ then we must have $y \in S$. Then either $r(y') = r(y)$, or,
                    since $X$ is convex, it follows from Remark~\ref{segSlope} 
                    that $r(y') \adj_{c_2} y = r(y)$.
          \end{itemize}
    \item Suppose $y$ is vertically above or below a point $x \in S$, so $r(y)=x$.
          Since $X$ is convex, it follows from Remark~\ref{segSlope} that $r(y') \adjeq_{c_2} x =r(y)$.
    \item Suppose $p_1(y) < m$. Since $X$ is convex, it follows from Remark~\ref{segSlope} 
          that $r(y') \adjeq_{c_2} r(y)$.
    \item Suppose $p_1(y) > M$. Since $X$ is convex, it follows from Remark~\ref{segSlope} 
          that $r(y') \adjeq_{c_2} r(y)$.
\end{itemize}
Thus $r \in C(X,c_2)$. Therefore, $r$ is a retraction. Clearly,
$r(Y \setminus Int(S)) = S$. This completes the proof.
\end{proof}

\begin{thm}
\label{diskAFPP}
Let $X \subset Y=[a,b]_{\Z} \times [c,d]_{\Z} \subset \Z^2$, such that $X$ is
digitally convex. Then $(X,c_2)$ has the AFPP.
\end{thm}

\begin{proof}
By Theorem~\ref{haveAFPP}(5), $(Y,c_2)$ has the AFPP. By Theorem~\ref{cvxRetractOfRect}, 
$(X,c_2)$ is a retract of $(Y,c_2)$. By Theorem~\ref{haveAFPP}(3), $(X,c_2)$ has the AFPP.
\end{proof}

\begin{cor}
Let $X = X_1 \times X_2$, where $X_1 \subset \Z^n$, $(X_1,c_n)$ has the AFPP,
$X_2 \subset \Z^2$, and $X_2$ is a digitally convex disk. Then
$(X,c_{n+2})$ has the AFPP.
\end{cor}

\begin{proof}
Clearly there exists $Y=[a,b]_{\Z} \times [c,d]_{\Z}$ such that $X_2 \subset Y$.
By Theorem~\ref{haveAFPP}(4), $(X_1 \times Y, c_{n+2})$ has the AFPP.
By Theorem~\ref{cvxRetractOfRect}, there is a $c_2$-retraction
$r: Y \to X_2$. Then
$\id_{X_1} \times r: X_1 \times Y \to X$ is a $c_{n+2}$-retraction.
The assertion follows from Theorem~\ref{haveAFPP}(3).
\end{proof}

\begin{prop}
\label{retractToBdHole}
Let $X \subset \Z^2$ be finite. Suppose $X' \subset \Z^2$ is a convex disk with bounding
curve~$S$ such that $X' \setminus Int(S)$
is a $c_2$-component of $\Z^2 \setminus X$. Then there is a
$c_2$-retraction of $X$ onto $S$.
\end{prop}

\begin{proof}
By Theorem~\ref{cvxRetractOfRect}, there is a $c_2$-retraction
$r: X \cup X' \to X'$ such that $r(X) = S$. Then
$r|_X: X \to S$ is a retraction.
\end{proof}

\begin{thm}
\label{nonAFPP}
Let $X \subset \Z^2$ be finite. Suppose $X' \subset X$ is a convex disk with bounding curve~$S$
such that $X' \setminus S$ is a 
$c_2$-component of $\Z^2 \setminus X$. Suppose there is a continuous
$F: S \to S$ such that  
\begin{equation}
\label{far}
\mbox{for each } x \in S, ~N^*(X,\kappa,x) \cap N^*(X,\kappa, F(x)) = \emptyset.
\end{equation}
Then $(X,c_2)$ does not have the AFPP.
\end{thm}

\begin{proof}
By Proposition~\ref{retractToBdHole}, there is a $c_2$-retraction $r: X \to S$. Let $F$ be as
described above and let $f: X \to X$ be the function $f(x)=F \circ r(x)$. Since composition
preserves continuity, we have $f \in C(X,c_2)$.

Consider the following cases.
\begin{itemize}
    \item Suppose $y \not \adjeq_{c_2} x$ for all $x \in S$. Then in particular,
          $y \not \adjeq_{c_2} f(y)$, so $y$ is not an approximate fixed point for $f$.
    \item Suppose $y \adjeq_{c_2} x$ for some $x \in S$. Then the continuity of $f$
          implies $f(y) \adjeq_{c_2} f(x)$. It follows from~(\ref{far}) that $y$
          is not an approximate fixed point for $f$.
\end{itemize}
Thus $f$ does not have an approximate fixed point. The assertion follows.
\end{proof}

\begin{figure}
    \centering
    \includegraphics[height=3in]{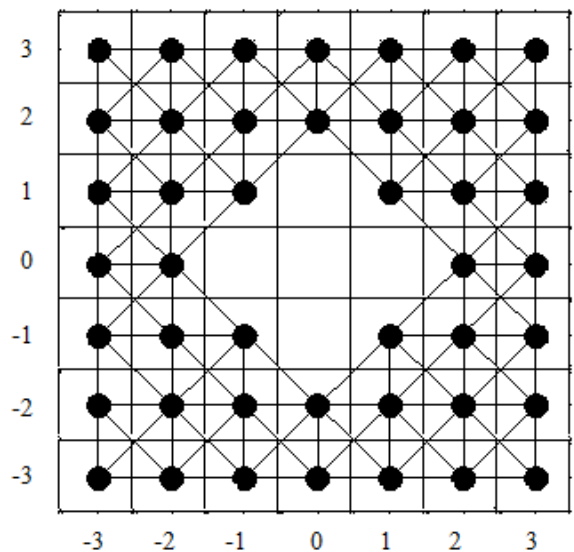}
    \caption{The image $X$ of Example~\ref{annulusExl}
       }
    \label{fig:annulus2}
\end{figure}

\begin{exl}
\label{annulusExl}
Let $X = [-3,3]_{\Z}^2 \setminus A$, where
$A = \{(x,y) \in \Z^2 ~ | ~ |x| + |y| \le 1\}$.
See Figure~\ref{fig:annulus2}.
As a bounding curve for $A$, we can take
$S = \{(x,y) \in \Z^2 ~ | ~ |x| + |y| = 2\}$.
Then $S$ is a $c_2$-simple closed curve. Let $F: S \to S$ be the map
$F(x,y) = (-x,-y)$. Then we may apply Theorem~\ref{nonAFPP} to conclude that
$(X,c_2)$ does not have the AFPP.
\end{exl}

\section{Further remarks}
We have explored relationships between the convexity of digital images in~$\Z^2$ and the AFPP.

In classical topology, every {\em absolute retract} (a contractible compactum with 
certain ``nice" local properties for which we need not consider analogs in
digital topology) has the FPP~\cite{Borsuk}. 
Since under the definition of digital homotopy in~\cite{Boxer99},
a digital simple closed curve of 4 points is digitally contractible~\cite{Boxer05},
Example~\ref{SCC4notAFPP} shows that contractibility based on~\cite{Boxer99} is not
sufficient for the AFPP. Recent papers of Staecker~\cite{StaeckerStrong} and Lupton, Oprea, and
Scoville~\cite{LuptonEtal} have developed a different notion of homotopy under which a
digital simple closed curve of 4 points is not digitally contractible; Staecker calls this
{\em strong homotopy}. This suggests the following questions concerning possible
extensions of Theorem~\ref{diskAFPP}.

\begin{question}
\label{contractibleAFPPq}
Let $X \subset \Z^2$ be finite and $c_2$-strongly contractible, i.e., contractible with
respect to strong homotopy. Does $(X,c_2)$ have the AFPP?
\end{question}

If Question~\ref{contractibleAFPPq} and the following Question~\ref{diskAFPPq} both have
affirmative answers, the latter result would be contained in the former.

\begin{question}
\label{diskAFPPq}
Let $X \subset \Z^2$ be a digital disk. Does $(X,c_2)$ have the AFPP?
\end{question}

\end{document}